\documentclass[12TP,draft]{article}

\begin{document}

\begin{center}
\LARGE\noindent\textbf{A theorem on even pancyclic bipartite digraphs}\\

\end{center}
\begin{center}
\noindent\textbf{Samvel Kh. Darbinyan }\\

Institute for Informatics and Automation Problems, Armenian National Academy of Sciences

E-mails: samdarbin@ipia.sci.am\\
\end{center}

\textbf{Abstract}

We prove that a strongly connected balanced bipartite directed graph of order $2a\geq 6$ with partite sets $X$ and $Y$ 
contains cycles of every length $2, 4, \ldots , 2a$, provided $d(x)+d(y)\geq 3a$ for every pair of vertices $x$, $y$ either both in $X$ or both in $Y$.

 
\section {Introduction} 

We consider directed graphs (digraphs) in the sense of \cite{[4]}. A cycle  is called Hamiltonian  if it  includes all the vertices of $D$. A digraph $D$ is Hamiltonian  if it contains a  Hamiltonian cycle.  
 There are numerous sufficient conditions for the existence of a Hamiltonian cycle in a digraph (see \cite{[2]} - \cite{[5]}, \cite{[7]}, \cite{[13]}, \cite{[14]}, \cite{[16]}, \cite{[17]}). The more general and classical ones are the following theorem by M. Meyniel.\\
 
\textbf{Theorem 1.1} (Meyniel \cite{[17]}). {\it Let $D$ be a strongly connected digraph of order $n\geq 2$. If $d(x)+d(y)\geq 2n-1$ for all pairs of non-adjacent vertices in $D$, then $D$ is Hamiltonian.}\\

Notice that Meyniel's theorem is a  generalization of Ghouila-Houri's and Woodall's theorems. 

A digraph is pancyclic if it contains  cycles of every length $k$,  $3\leq k\leq n$,
 where $n$ is the order of $D$.  
There are various sufficient conditions for a digraph  (undirected graph) to be Hamiltonian are also  sufficient for the digraph to be pancyclic (see \cite{[4]}, \cite{[7]}, \cite{[9]},  \cite{[10]}, \cite{[14]}, \cite{[15]}, \cite{[16]}, \cite{[19]}). In  \cite{[8]} and  \cite{[9]}, the author studded the pancyclcity of a digraph with the condition of the Meyniel
 theorem. Before stating the main result of \cite{[9]} we need define a family of digraphs.\\

\textbf{Definition 1.2}. {\it For any integers $n$ and $m$, $(n+1)/2 < m\leq n-1$, let $\Phi_n^m$ denote the set of  
 digraphs $D$, which satisfy the following conditions: (i) $V(D)= \{x_1,x_2,\ldots , x_n\}$; (ii) $x_nx_{n-1}\ldots x_2x_1x_n$ is a Hamiltonian cycle in $D$; (iii) for each $k$, $1\leq k\leq n-m+1$, the vertices $x_k$ and $x_{k+m-1}$ are not adjacent; 
(iv) $x_jx_i\notin A(D)$ whenever $2\leq i+1< j\leq n$ and (v) the sum of degrees for any two distinct non-adjacent vertices at least $2n-1$.}\\

\textbf{Theorem 1.3} (Darbinyan \cite{[9]}). {\it Let $D$ be a strongly connected digraph of order $n\geq 2$. Suppose that $d(x)+d(y)\geq 2n-1$ for all pairs of distinct non-adjacent vertices $x$, $y$ in $D$. Then either (a) $D$ is pancyclic or (b) $n$ is even and $D$ is isomorphic to one of $K^*_{n/2,n/2}$, $K^*_{n/2,n/2}\setminus \{e\}$, where $e$ is an arbitrary arc of $K^*_{n/2,n/2}$, or (c) $D\in \Phi^m_n $ (in this case $D$ does not contain a cycle of length $m$).}\\

Theorem 1.3, later also was proved independently by Benhocine \cite{[6]}. 
In \cite{[5]}, Bang-Jensen, Gutin and Li conjectured the following strengthening of Meyniel's theorem.\\

\textbf{Conjecture 1.4}. {\it Let $D$ be a strongly connected digraph of order $n$. Suppose that $d(x)+d(y)\geq 2n-1$ for every pair of non-adjacent distinct vertices $x$, $y$ with a common out-neighbor or a common in-neighbor. Then $D$ is Hamiltonian.}\\

They also conjectured that this can even be generalized to\\

\textbf{Conjecture 1.5}. {\it Let $D$ be a strongly connected digraph of order $n$. Suppose that $d(x)+d(y)\geq 2n-1$ for every pair of non-adjacent distinct vertices $x$, $y$ with a common  in-neighbor. Then $D$ is Hamiltonian.}\\

In \cite{[3]} and \cite{[5]}, it was proved that   Conjecture 1.4 (1.5) is true if we also require an additional condition.\\

\textbf{Theorem 1.6} (Bang-Jensen, Guo, Yeo \cite{[3]}). {\it Let $D$ be a strongly connected digraph of order $n\geq 2$. Suppose that $min\{d^+(x)+d^-(y), d^-(x)+d^+(y) \}\geq n-1$ and  $d(x)+d(y)\geq 2n-1$ for any pair of non-adjacent vertices $x,y$ with a common in-neighbor or a common out-neighbor. Then $D$ is Hamiltonian.}\\

\textbf{Theorem 1.7} (Bang-Jensen, Gutin, Li \cite{[5]}). {\it Let $D$ be a strongly connected digraph of order $n\geq 2$. Suppose that $min\{d(x),d(y)\}\geq n-1$ and  $d(x)+d(y)\geq 2n-1$ for any pair of non-adjacent vertices $x,y$ with a common in-neighbor. Then $D$ is Hamiltonian.}\\

In \cite{[3]},  also it was proved that if in Conjecture 1.4  we replace $5n/2-4$ instead of $2n-1$, then  Conjecture 1.4 is true.

A digraph $D$ is a bipartite  if there exists a partition $X$, $Y$ of its vertex set into two partite sets such that every arc of 
$D$ has its end-vertices in different partite sets. 
It is called balanced if $|X|=|Y|$.

An  analogue of Meyniel's theorem for the hamiltonicity of balanced digraphs was given by   Adamus, Adamus and Yeo \cite{[2]}. \\

\textbf{Theorem 1.8} (Adamus, Adamus, Yeo \cite{[2]}). {\it Let $D$ be a balanced bipartite digraph of order $2a\geq 4$ with partite sets $X$ and $Y$. Then $D$ is Hamiltonian provided one of the following holds:
 
(a) $d(x)+d(y)\geq 3a+1$ for each pair of non-adjacent vertices $x,y\in X\cup Y$; 

(b) $D$ is strongly connected and $d(x)+d(y)\geq 3a$ for each pair of non-adjacent vertices $x,y\in X\cup Y$;

(c) the minimal degree of $D$ at least $(3a+1)/2$;

(d) $D$ is strongly connected and its the minimal degree  at least $3a/2$.}\\

Meszka  \cite{[18]} proved  the following theorem.

\textbf{Theorem 1.9} (Meszka  \cite{[18]}). {\it Let $D$ be a balanced bipartite digraph of order $2a\geq 4$ with partite sets 
$X$ and $Y$. Suppose that  $d(x)+d(y)\geq 3a+1$ for each two  vertices $x,y$ either both in $X$ or both in $Y$. Then $D$  
 contains cycles of all even lengths $4, 6, \ldots , 2a$.}\\

Notice that Theorem 1.8(a) is an immediate corollary of Theorem 1.9.
There are some versions of Conjecture 1.4 and  1.5 for balanced bipartite digraphs. (Theorems 1.10-1.11 below).\\

\textbf{Theorem 1.10} (Adamus \cite{[1]}). {\it Let $D$ be a strongly connected balanced bipartite digraph of order $2a\geq 6$. If $d(x)+d(y)\geq 3a$ for every pair of distinct vertices $x, y$ with a common out-neighbor or a common in-neighbor, then $D$ is Hamiltonian.}\\   

 The last theorem improved Theorem 1.8. 
An analogue of Theorem 1.6   was given by  Wang \cite{[20]}, and recently strengthened by the author \cite{[11]}.\\

\textbf{Theorem 1.11} (Wang \cite{[20]}). {\it Let $D$ be a strongly connected balanced bipartite digraph of order $2a\geq 4$. Suppose that, for every dominating pair of vertices $\{x,y\}$, either $d(x)\geq 2a-1$ and $d(y)\geq a+1$ or $d(y)\geq 2a-1$ and $d(x)\geq a+1$. Then $D$ is Hamiltonian}.\\

Before stating the next theorems we need to define a digraph of order eight.\\

\textbf{Definition 1.12.} {\it Let $D(8)$ be the  bipartite digraph  with partite sets $X=\{x_0,x_1,x_2,x_3\}$ and 
$Y=\{y_0,y_1,y_2,y_3\}$, and  $A(D(8))$ contains exactly the arcs $y_0x_1$, $y_1x_0$, $x_2y_3$, $x_3y_2$ and all the arcs of the following 2-cycles: 
$x_i\leftrightarrow y_i$, $i\in [0,3]$, $y_0\leftrightarrow x_2$, $y_0\leftrightarrow x_3$, $y_1\leftrightarrow x_2$ and  $y_1\leftrightarrow x_3$.}\\

It is not difficult to check that $D(8)$ is strongly connected, $max \{d(x), d(y)\}\geq 2a-1$ for every  pair of vertices $\{x,y\}$ with a common out-neighbor and it is not Hamiltonian. (Indeed, if $C$ is a Hamiltonian cycle in $D(8)$, then $C$ would  contain the arcs $x_1y_1$ and $x_0y_0$ and therefore, the path $x_1y_1x_0y_0$ or the path $x_0y_0x_1y_1$, which is impossible since $N^-(x_0)=N^-(x_1)=\{y_0,y_1\}$).\\

\textbf{Theorem 1.13} (Darbinyan \cite{[11]} and \cite{[12]}). {\it Let $D$ be a strongly connected balanced bipartite digraph of order $2a\geq 8$ with partite sets $X$ and $Y$.  If $D$ is not a directed cycle and  $ max\{d(x), d(y)\}\geq 2a-1$ for every  pair of distinct vertices $\{x,y\}$ with a common out-neighbor,
  then either $D$ contains  cycles of all even lengths less than or equal to $2a$ or $D$ is isomorphic to the digraph $D(8)$.}\\

Notice that Theorem 1.11 is an immediate consequence of Theorem 1.13.
 
In this paper using Theorem 1.3 and some arguments of \cite{[18]} we prove the following theorem.\\

\textbf{Theorem 1.14.} {\it Let $D$ be a strongly connected balanced bipartite digraph of order $2a\geq 6$ with partite sets $X$ and $Y$.  If  $d(x)+ d(y)\geq 3a$ for every  pair of distinct vertices $\{x,y\}$ either both in $X$ or both in 
$Y$, 
  then $D$ contains  cycles of all even lengths less than or equal to $2a$. }\\
 
\textbf{Remark}. {\it For any integer $a\geq 2$ there is a non-strongly connected balanced bipartite digraph $D$ of order $2a$ with partite sets $X$ and $Y$, such that  $d(x)+ d(y)\geq 3a$ for every  pair of distinct vertices $\{x,y\}$ either both in $X$ or both in 
$Y$. }

To see this, we take two   balanced bipartite complete digraph  of order $a$ ($a$ is even) with partite sets $U$, $V$ and $Z$, $W$, respectively. By adding all the possible arcs from $Z$ to $V$ and  from $W$ to $U$ we obtain a digraph $D$. It is easy to check that 
$d(x)+ d(y)\geq 3a$ for every  pair of non-adjacent distinct vertices $\{x,y\}$ of $D$,
 but $D$ is not strongly connected. \\

\section {Terminology and Notation}

  In this paper we consider finite digraphs without loops and multiple arcs. Terminology  and notation not described below follow \cite{[4]}. 
The vertex set and the arc set of a digraph $D$ are    denoted 
  by $V(D)$  and   $A(D)$, respectively.  The order of $D$ is the number
  of its vertices. 
 If $xy\in A(D)$, then we also write $x\rightarrow y$ and say that $x$ dominates $y$ or $y$ is an out-neighbor of $x$ and $x$ is an in-neighbor of $y$. If $x\rightarrow y$ and $y\rightarrow x$ we shall use the notation $x\leftrightarrow y$  ($x\leftrightarrow y$ is called  2-cycle). Let $\overrightarrow {a} [x,y]=1$ if $xy\in A(D)$ and $\overrightarrow{a} [x,y]=0$ if $xy\notin A(D)$.

 If $A$ and $B$ are two disjoint subsets of $V(D)$ such that every
   vertex of $A$ dominates every vertex of $B$, then we say that $A$ dominates $B$, denoted by $A\rightarrow B$. 
 Similarly, $A\leftrightarrow B$ means that 
$A\rightarrow B$ and $B\rightarrow A$. If $x\in V(D)$
   and $A=\{x\}$ we sometimes write $x$ instead of $\{x\}$.
Let $N_D^+(x)$, $N_D^-(x)$ denote the set of  out-neighbors, respectively the set  of in-neighbors of a vertex $x$ in a digraph $D$.  If $A\subseteq V(D)$, then $N_D^+(x,A)= A \cap N_D^+(x)$ and $N_D^-(x,A)=A\cap N_D^-(x)$. 
The out-degree of $x$ is $d_D^+(x)=|N_D^+(x)|$ and $d_D^-(x)=|N_D^-(x)|$ is the in-degree of $x$. Similarly, $d_D^+(x,A)=|N_D^+(x,A)|$ and $d_D^-(x,A)=|N_D^-(x,A)|$. The degree of the vertex $x$ in $D$ is defined as $d_D(x)=d_D^+(x)+d_D^-(x)$ (similarly, $d_D(x,A)=d_D^+(x,A)+d_D^-(x,A)$). We omit the subscript if the digraph is clear from the context. The subdigraph of $D$ induced by a subset $A$ of $V(D)$ is denoted by $D\langle A\rangle$ or $\langle A\rangle$ for brevity. 

For integers $a$ and $b$, $a\leq b$, let $[a,b]$  denote the set of
all the integers which are not less than $a$ and are not greater than
$b$.

 The path (respectively, the cycle) consisting of the distinct vertices $x_1,x_2,\ldots ,x_m$ ( $m\geq 2 $) and the arcs $x_ix_{i+1}$, $i\in [1,m-1]$  (respectively, $x_ix_{i+1}$, $i\in [1,m-1]$, and $x_mx_1$), is denoted by  $x_1x_2\cdots x_m$ (respectively, $x_1x_2\cdots x_mx_1$). The length of a cycle or a path is the number of its arcs.
We say that $x_1x_2\cdots x_m$ is a path from $x_1$ to $x_m$ or is an $(x_1,x_m)$-path. If a digraph $D$  contains a path from a vertex $x$ to a vertex $y$ we say that $y$ is reachable from  $x$ in $D$. In particular, $x$ is reachable from  itself.

Let $K^*_{a,b}$ denote the complete bipartite digraph with partite sets of cardinalities $a$ and $b$.

A digraph $D$ is strongly connected (or, just, strong) if there exists a path from $x$ to $y$ and a path from $y$ to $x$ for every pair of distinct vertices $x,y$.
    
Two distinct vertices $x$ and $y$ are adjacent if $xy\in A(D)$ or $yx\in A(D) $ (or both).

Let $D$ be a  bipartite digraph with  partite sets $X$ and $Y$. A matching from $X$ to $Y$ (from $Y$ to $X$) is an independent set of arcs with origin in $X$ and terminus in $Y$ (origin in $Y$ and terminus in $X$). (A set of arcs with no common end-vertex is called  independent). If $D$ is balanced, one says that such a matching is perfect if it consists of precisely $|X|$ arcs. 

\section { Preliminaries }

Let us recall some results (Lemmas 3.1 and 3.4) which will be used in this paper.\\

  \textbf{Lemma 3.1} (Adamus \cite{[1]}). {\it Let $D$ be a strongly connected balanced bipartite digraph of order $2a\geq 4$. If 
$d(x)+d(y)\geq 3a$ for every pair of vertices $x$, $y$ with a common in-neighbor or a common out-neighbor, then  $D$ contains a perfect matching from $Y$ to $X$ and a perfect matching from $X$ to $Y$.} \\

\textbf{Definition 3.2.}  {\it Let $D$ be a  balanced bipartite digraph of order $2a\geq 6$  with partite sets $X$ and $Y$. Let $k\geq 0$ be an integer. We will say that $D$  satisfies condition $A_k$ when $d(x)+d(y)\geq  3a+k$ for every pair of distinct  vertices $x$, $y$ either both in $X$ or both in $Y$.}\\

Follows \cite{[18]}, we give the following definition.
 
\textbf{Definition 3.3.} {\it Let $D$ be a  balanced bipartite digraph of order $2a\geq 6$  with partite sets $X$ and $Y$. Let $M_{y,x}= \{y_ix_i\in A(D)\, |\, i=1, 2, \ldots , a\}$ be an arbitrary perfect matching from $Y$ to $X$. We define a digraph
$D^*[M_{y,x}]$ of order $a$ and with vertex set $\{v_1,v_2, \ldots , v_a\}$ as follows: Each vertex $v_i$ corresponds to a pair $\{x_i,y_i\}$ of vertices in $D$ and for each pair of distinct vertices $v_l, v_j$,  $v_lv_j\in A(D^*[M_{y,x}])$ if and only if $x_ly_j\in A(D)$.}
\\\

 Let $D$ be a  balanced bipartite digraph  with partite sets $X$ and $Y$. Let $M_{y,x}$ be a perfect matching from $Y$ to $X$ in $D$ and $D^*[M_{y,x}]$ be its corresponding digraph. Further, in this paper we will denote the vertices of $D$ (respectively, of $D^*[M_{y,x}]$) by letters $x$, $y$ (respectively, $u$, $v$) with subscripts or without subscripts.

 The size of a perfect matching $M_{y,x}= \{y_ix_i\in A(D)\, |\, i=1, 2, \ldots , a\}$ from $Y$ to $X$ in $D$ (denoted by $s(M_{y,x})$) is the number of arcs $y_ix_i$ such that $x_iy_i\notin A(D)$.\\

Using \cite{[18]}, we can formulate the following lemma.\\

\textbf{Lemma 3.4} (Meszka \cite{[18]}). {\it Let $D$ be a  balanced bipartite digraph of order $2a\geq 6$  with partite sets $X$ and $Y$. Let $M_{y,x}= \{y_ix_i\in A(D)\, |\, i=1, 2, \ldots , a\}$ be an arbitrary perfect matching from $Y$ to $X$. Then the following holds.

(i). $d^+(v_i)= d^+(x_i)-\overrightarrow{a}[x_i,y_i]$     and $d^-(v_i)= d^-(y_i)-\overrightarrow{a}[x_i,y_i]$.

(ii). If  $D^*[M_{y,x}]$ contains a cycle of length $k$, where $k\in [2,a]$, then $D$ contains a cycle of length $2k$.

(iii). Let $a\geq 4$ is even and $D^*[M_{y,x}]$  is isomorphic to $K^*_{a/2,a/2}$
 with partite sets $\{v_1,v_2,\ldots , v_{a/2}\}$ and
 $\{v_{a/2+1},v_{a/2+2},\ldots , v_{a}\}$. If $D$ contains an arc from $\{y_1,y_2,\ldots , y_{a/2}\}$ to  $\{x_{a/2+1},x_{a/2+2},\ldots , x_{a}\}$, say $y_{a/2} x_a\in A(D)$. 
Then $D$ contains  cycles of every length $2k$, $k=2,3,  \ldots, a$}.

\textbf{Proof of Lemma 3.4.} The proof of Lemma 3.4 can be found in \cite{[18]}, but we give it here for completeness.

(i). It follows immediately from the definition of $D^*[M_{y,x}]$.

(ii). Indeed, if $v_{i_1}v_{i_2}\ldots v_{i_k}v_{i_1}$ is a cycle of length $k$ in $D^*[M_{y,x}]$, then
 $y_{i_1}x_{i_1}y_{i_2}x_{i_2}y_{i_3}\ldots y_{i_k}x_{i_k}y_{i_1}$ is a cycle of length $2k$ in $D$.

(iii).  By (ii), it is clear that $D$ contains cycles of every length $4k$, $k=1,2 \ldots , a/2$. It remains to show that $D$ also contains cycles of every length $4k+2$, $k=1,2 \ldots , a/2-1$.  Indeed, since $x_iy_j \in A(D)$ and $x_jy_i\in A(D)$ for all
 $ i\in [1, a/2]$, $j\in [a/2+1,a]$ and  $y_{a/2} x_a\in A(D)$, from the definition of $D^*[M_{y,x}]$  
 it follows that $ y_1x_1y_{a/2+1}x_{a/2+1} y_2 x_2y_{a/2+2}x_{a/2+2}y_3x_3 \ldots 
 $ $x_ky_{a/2+k}x_{a/2+k}y_{a/2}x_{a}y_1$ is a cycle of length $4k+2$ in $D$.                     \fbox \\\\

\section {Proof of the main result}

The  proof of Theorem 1.14 will be based on the  following three lemmas.\\

\textbf{Lemma 4.1.} {\it Let $D$ be a strongly connected balanced bipartite digraph of order $2a\geq 6$  with partite sets $X$ and $Y$.  If $D$  satisfies condition $A_0$, then $D$ contains  cycles of  lengths 2 and 4}. 

\textbf{Proof of Lemma 4.1}. From condition $A_0$ immediately follows that $D$ contains a cycle of length 2. We will prove that $D$ contains a cycle of length 4. By Lemma 3.1, $D$ contains a perfect matching from $Y$ to $X$. Let $M_{y,x}= \{y_ix_i\in A(D)\, |\, i=1, 2, \ldots , a\}$ be an arbitrary perfect matching from $Y$ to $X$. If for some 
integers $i,j$, $1\leq i\not= j\leq a$, the arcs $x_iy_j$, $x_jy_i$ are in $D$, then $x_iy_jx_jy_ix_i$ is a cycle of length 4. We may therefore assume that
for any pair of integers $i,j$, $1\leq i\not= j\leq a$,
$$
\overrightarrow{a}[x_i,y_j]+\overrightarrow{a}[x_j,y_i]\leq 1.   
$$
Therefore, for all $i\in [1,a]$ we have
$$ 
d^-(y_i)\leq a-d^+(x_i)-1, \; \hbox{if}\;\overrightarrow{a}[x_i,y_i]=0, \;\hbox{and} \;d^-(y_i)\leq a-d^+(x_i)+1, \; \hbox{if}\;  \overrightarrow{a}[x_i,y_i]=1.  \eqno (1)
$$
Assume that there are two distinct integers $i, j$, $1\leq i,j\leq a$, such that 
$\overrightarrow{a}[x_i,y_i]=\overrightarrow{a}[x_j,y_j]=0$. Then, by (1), 
$d^-(y_i)+d^+(x_i)\leq a-1$ and $d^-(y_j)+d^+(x_j)\leq a-1$. These together with condition $A_0$ and the fact that the semi-degrees of every vertex in $D$ are bounded above by $a$ thus imply that 
$$
6a\leq d(x_i)+d(x_j)+d(y_i)+d(y_j)=d^-(y_i)+d^+(x_i)+d^-(y_j)+d^+(x_j)$$ $$+ d^+(y_i)+d^+(y_j)+d^-(x_i)+d^-(x_j)\leq 6a-2,
$$
which is a contradiction.

Assume now that for some $i\in [1,a]$, $\overrightarrow{a}[x_i,y_i]=0$ and for all $j\in [1,a]\setminus \{i\}$,
 $\overrightarrow{a}[x_j,y_j]= 1$. Without loose of generality we may assume that $i=1$ and $j=2$. By (1),   $d^-(y_1)+d^+(x_1)\leq a-1$ and $d^-(y_2)+d^+(x_2)\leq a+1$. If for some $k\in [3,a]$, $y_2x_k\in A(D)$ and $y_kx_2\in A(D)$, then $x_2y_2x_ky_kx_2$ is a cycle of length 4 in $D$. We may therefore assume that  
$\overrightarrow{a}[y_2,x_k]+\overrightarrow{a}[y_k,x_2]\leq 1$ for all $k\in [3,a]$. This implies that 
$$
d^-(x_2)\leq a-2-d^+(y_2, X\setminus \{x_1,x_2\})+2,
$$ 
and hence $d^-(x_2)+d^+(y_2)\leq a+2$. Using the above inequalities and condition $A_0$, we obtain
$$
6a\leq d(x_1)+d(x_2)+d(y_1)+d(y_2)=d^-(y_1)+d^+(x_1)+d^-(y_2)+d^+(x_2)$$ $$+ d^-(x_2)+d^+(y_2)+d^-(x_1)+d^+(y_1)\leq 5a+2,
$$
which is a contradiction since $a\geq 3$.

Assume finally that $x_iy_i\in A(D)$ for all $i\in [1,a]$. In this case, by the symmetry between the vertices $x_i$ and $y_i$, similarly to (1), we obtain that $d^-(x_i)+d^+(y_i)\leq a+1$. This together with (1) implies that for any $i$, $j$ ($1\leq i\not= j\leq a$), 
$$
6a\leq d(x_i)+d(x_j)+d(y_i)+d(y_j)\leq 4a+4,
$$ a contradiction since $a\geq 3$. Lemma 4.1 is proved. \fbox \\\\

\textbf{Lemma 4.2}. {\it Let $D$ be a strongly connected balanced bipartite digraph of order $2a\geq 6$  with partite sets $X$ and $Y$. 
Let $M_{y,x}= \{y_ix_i\in A(D)\, |\, i=1, 2, \ldots , a\}$ be a perfect matching from $Y$ to $X$ in $D$ such that the size of $M_{y,x}$ is maximum among the sizes of all the perfect matching from $Y$ to $X$ in $D$. 
 If $D$  satisfies condition $A_0$, then the digraph 
$D^*[M_{y,x}]$ 
is strongly connected or $D$ contains cycles of all lengths $2, 4, \ldots , 2a$}. 

\textbf{Proof of Lemma 4.2}. Notice that, by Lemma 4.1, $D$ contains cycles of lengths 2 and 4. Suppose that the digraph $D^*[M_{y,x}]$ 
is not strongly connected.  Then in $D^*[M_{y,x}]$ there are two distinct vertices, say $v_1$ and $v_j$, such that there is no path from $v_1$ to $v_j$ in $D^*[M_{y,x}]$. Let  $U$ be the set of all vertices reachable from $v_1$ and $W$ be the set of all vertices from which $v_j$ is reachable.  Notice that $v_1\in U$, $v_j\in W$ and $U\cap W= \emptyset$. 

\textbf{Case 1}. $d^+(v_1)\geq 1$ and $d^-(v_j)\geq 1$.

Then $|U|\geq 2$ and $|W|\geq 2$. Let $v_l$, $v_k$ be arbitrary two distinct vertices in $U$ and   $v_p$, $v_q$  arbitrary two distinct vertices in $W$. From condition $A_0$ and the fact that the semi-degrees of every vertex in $D$ are bounded above by $a$ it follows that
$$
d^+(x_l)+d^+(x_k)\geq a \quad \hbox{and} \quad d^-(x_p)+d^-(x_q)\geq a.  \eqno (2)
$$
By Lemma 3.2(i), 
$$
 d^+(v_l)+d^+(v_k)=d^+(x_l)+d^+(x_k)- \overrightarrow{a}[x_l,y_l]- \overrightarrow{a}[x_k,y_k],
$$
and
$$
 d^-(v_p)+d^-(v_q)= d^-(y_p)+d^-(y_q)- \overrightarrow{a}[x_p,y_p]- \overrightarrow{a}[x_q,y_q]. \eqno (3)
$$
These and (2) imply that $d^+(v_l)+d^+(v_k)\geq a-2$ and $d^-(v_p)+d^-(v_q)\geq a-2$, which in turn implies that 
$|U|\geq a/2$ and $|W|\geq a/2$.

If  $d^+(v_l)+d^+(v_k)\geq a-1$ or $d^-(v_p)+d^-(v_q)\geq a-1$, then $|U|\geq (a+1)/2$ or $|W|\geq (a+1)/2$, respectively. Hence
$|U|+|W|\geq (2a+1)/2$, which is a contradiction since $|U|+|W|\leq a$. Using (2) and (3), we may therefore assume that
$$
 d^+(v_l)+d^+(v_k)= d^+(x_l)+d^+(x_k)-2= d^-(v_p)+d^-(v_q)= d^-(y_p)+d^-(y_q)-2=a-2.
$$
Then it is easy to see that the arcs $x_ly_l$, $x_ky_k$, $x_py_p$ and $x_qy_q$ are in $D$, $|U|=|W|=a/2$ and 
$V(D^*[M_{y,x}])= U\cup W$. In particular, $a$ is even. Without  loss of generality we assume that $U=\{v_1,v_2, \ldots , v_{a/2}\}$ and 
$W=\{v_{a/2+1},v_{a/2+2}, \ldots , v_{a}\}$. Since there is no arc from a vertex in $U$ to a vertex in $W$, the following holds:
$$
A(\{x_{1},x_{2}, \ldots , x_{a/2}\}\rightarrow \{y_{a/2+1},y_{a/2+2}, \ldots , y_{a}\})=\emptyset. \eqno (4)
$$
Therefore, if $i\in [1,a/2]$ and $j\in [a/2+1,a]$, then $d^+(x_i)\leq a/2$ and $d^-(y_j)\leq a/2$. 
These together with (2)
 imply that $d^+(x_i)=d^-(y_j)= a/2$ and 
$$x_i\rightarrow \{y_{1},y_{2}, \ldots , y_{a/2}\} \quad \hbox{ and} \quad  
\{x_{a/2+1},x_{a/2+2}, \ldots , x_{a}\}\rightarrow y_j \eqno (5)
$$
for all  $i\in [1,a/2]$ and $j\in [a/2+1,a]$, respectively. Therefore, by condition $A_0$, 
$$
3a\leq d(x_i)+d(x_k)\leq a+d^-(x_i)+d^-(x_k),
$$
 for any pair of $i,k\in [1,a/2]$. This  implies that $d^-(x_i)=d^-(x_k)=a$, i.e.,
$\{y_1,y_2,\ldots , y_a\}\rightarrow \{x_i,x_k\}$. Similarly,  
$y_j\rightarrow \{x_{1},x_2\ldots ,x_a\}$, for all $j\in [a/2+1,a]$. From this and  (5) it follows  that the induced subdigraphs 
$\langle\{x_1, x_2,\ldots ,x_{a/2},y_1, y_2,$ $\ldots ,y_{a/2}\}\rangle$ and 
$\langle\{x_{a/2+1},x_{a/2+2},\ldots ,x_a,y_{a/2+1},y_{a/2+2},\ldots ,y_a\}\rangle $ are balanced bipartite complete digraphs.
Therefore, $D$ contains cycles of all lengths $2,4,\ldots ,a$. It remains to show that $D$ also contains cycles of every length $a+2b$, $b\in [1,a/2]$.   Since $D$ is strongly connected and (4), it follows that there is an arc from a vertex in $\{y_1, y_2,\ldots ,y_{a/2}\}$ to a vertex in  $\{x_{a/2+1},x_{a/2+2},\ldots ,x_a\}$. Without loss of generality we may assume that $y_{a/2}x_{a/2+1}\in A(D)$. Then
 $x_1y_1x_2y_2\ldots x_{a/2}y_{a/2}x_{a/2+1}$ $y_{a/2+1}x_{a/2+2}\ldots x_{a/2+b}y_{a/2+b}x_{1}$ is a cycle of length $a+2b$.
 Thus $D$ contains cycles of all lengths $2,4,\ldots , 2a$. This completes the discussion of   Case 1.

\textbf{Case 2}. $d^+(v_1)=0$.
 
Then $d^+(x_1)=1$ and $x_1y_1\in A(D)$, since $D$ is strongly connected. Hence $d(x_1)\leq a+1$. This together with condition $A_0$ implies that $a\leq d(x_1)\leq a+1$. We now distinguish two subcases.

\textbf{Case 2.1}. $d(x_1)=a$.

Then $d(x_i)\geq 2a$ for all $i\in [2,a]$ because of condition $A_0$. Therefore, the induced subdigraph $\langle Y\cup X\setminus \{x_1\}\rangle$ is a complete bipartite digraph with partite sets $Y$ and $X\setminus \{x_1\}$. It is clear that $D$ contains cycles of every lengths $2, 4, \ldots , 2a-2$. Since $d(x_1)=a$, $d^+(x_1)=1$ and $a\geq 3$, we have that $d^-(x_1)=a-1\geq 2$. Without loose of generality we may assume that $y_2x_1\in A(D)$. Then
 $x_2y_2x_1y_1x_3y_3\ldots x_ay_ax_2$ is a cycle of length $2a$.

\textbf{Case 2.2}. $d(x_1)=a+1$.

Then $\{y_1,y_2,\ldots , y_a\}\rightarrow x_1$ because of $d^+(x_1)=1$, and, by condition $A_0$, $d(x_i)\geq 2a-1$ for all $i\in [2,a]$.  Observe that if for some 
$i\in [2,a]$,  $y_1x_i\in A(D)$, then 
$M^i_{y,x}:=\{y_ix_1,y_1x_i\}\cup \{y_jx_j \, |\, j\in [1,a]\setminus \{1,i\}\}$ is a perfect matching from $Y$ to $X$ in $D$.

Assume that for some $i\in [2,a]$,  $x_iy_1\notin A(D)$. Then $y_1x_i\in A(D)$ because of $d(x_i)\geq 2a-1$.    Since $x_1y_1\in A(D)$, $x_iy_1\notin A(D)$ and  $x_1y_i\notin A(D)$, it follows that  $s(M^i_{y,x})>s(M_{y,x})$, which contradicts the choice of $M_{y,x}$.
 We may therefore assume that $\{x_2,x_3,\ldots , x_a\}\rightarrow y_1$. If $y_1x_i\in A(D)$ and  
$x_iy_i\in A(D)$, where $i\in [2,a]$, then again we have $s(M^i_{y,x})>s(M_{y,x})$, since the arcs $x_1y_1$, $x_iy_i$ are in $D$ and $x_1y_i\notin A(D)$.  We may therefore assume that $\overrightarrow{a}[y_1,x_i]+\overrightarrow{a}[x_i,y_i]\leq 1$.
This together with $d(x_i)\geq 2a-1$, $i\in [2,a]$, thus imply that 
$$
 \{y_2,y_3,\ldots , y_a\}\rightarrow x_i\rightarrow \{y_2,y_3,\ldots , y_a\}\setminus \{y_i\}. \eqno (6)
$$
From strongly connectedness of $D$ and $d^+(x_1,\{y_2,y_3,\ldots , y_a\}=0$ it follows that $d^+(y_1, \{x_2,x_3, \ldots , x_a\})\not= 0$. Without loss of generality we assume that 
$y_1x_2\in A(D)$. Then, since $y_2x_1\in A(D)$ and (6), $x_1y_1x_2y_3x_3\ldots x_{k-1}y_kx_ky_2x_1$ is a cycle of length $2k$ for every $k\in [3,a]$. Lemma 4.2 is proved. \fbox \\\\

\textbf{Lemma 4.3.} {\it Let $D$ be a strongly connected balanced bipartite digraph of order $2a\geq 6$  with partite sets $X$ and $Y$.  Let $M_{y,x}= \{y_ix_i\in A(D)\, |\, i=1, 2, \ldots , a\}$ be a perfect matching from $Y$ to $X$ in $D$ such that the size of $M_{y,x}$ is maximum  among the sizes of all the perfect matching from $Y$ to $X$ in $D$.  
If $D$  satisfies condition $A_0$, then either $d(u)+d(v)\geq 2a-1$ for every pair of non-adjacent vertices $u$, $v$ in
$D^*[M_{y,x}]$ or $a\leq 4$ in which case $D$ contains cycles of all lengths $2, 4, \ldots , 2a$.} 

\textbf{Proof of Lemma 4.3}. Let $v_i$ and $v_j$ be two arbitrary distinct vertices in  $D^*[M_{y,x}]$. Put 
$$ g(i,j):=d^+(x_i)+d^+(x_j)+d^-(y_i)+d^-(y_j) \; \hbox{and} \; f(i,j):=d^-(x_i)+d^-(x_j)+d^+(y_i)+d^+(y_j).$$
By Lemma 3.2(i), we have  
$$
d(v_i)+d(v_j)=g(i,j)-2\overrightarrow{a}[x_i,y_i]-2\overrightarrow{a}[x_j,y_j]. \eqno (7)
$$
By condition $A_0$, we have
$$
6a\leq d(x_i)+d(x_j)+d(y_i)+d(y_j)=f(i,j)+g(i,j).          
$$
Hence $g(i,j)\geq 2a$, since the semi-degrees of every vertex of $D$ are bounded  above by $a$, and
$$
 f(i,j)\geq 6a-g(i,j). \eqno (8)
$$
Now we prove the following claim.

\textbf{Claim 1}. Assume that the vertices $v_i$ and $v_j$ in  $D^*[M_{y,x}]$ are not adjacent. Then the following holds:

(i). If $x_iy_i\in A(D)$ or $x_jy_j\in A(D)$, then $\overrightarrow{a}[y_i,x_j]+\overrightarrow{a}[y_j,x_i]\leq 1$.

(ii). If $x_iy_i\notin A(D)$ or $x_jy_j\notin A(D)$, then $d(v_i)+d(v_j)\geq 2a-1$ in  $D^*[M_{y,x}]$.

\textbf{Proof of the claim}. Since the vertices  $v_i$ and $v_j$ in  $D^*[M_{y,x}]$ are not adjacent, it follows that $x_iy_j\notin A(D)$ and $x_jy_i\notin A(D)$.

 (i). Suppose, to the contrary, that $x_iy_i\in A(D)$ or $x_jy_j\in A(D)$, but $\overrightarrow{a}[y_i,x_j]+\overrightarrow{a}[y_j,x_i]=2$. Then  
$M'_{y,x}=\{y_ix_j,y_jx_i\}\cup \{y_kx_k \, |\, k\in [1,a]\setminus \{i,j\}\}$ is a new perfect matching from $Y$ to $X$ in $D$. Since 
$x_jy_i\notin A(D)$, $x_iy_j\notin A(D)$ and $x_iy_i\in A(D)$ or $x_jy_j\in A(D)$, it follows that 
 $s(M'_{y,x})>s(M_{y,x})$, which contradicts the choice of $M_{y,x}$. 

(ii). If $\overrightarrow{a}[x_i,y_i]=\overrightarrow{a}[x_j,y_j]=0$, then from (7) and $g(i,j)\geq 2a$ it follows that
 $d(v_i)+d(v_j\geq 2a$. 
 We may therefore assume that $x_iy_i\in A(D)$. Then   $x_jy_j\notin A(D)$ by the assumption of Claim 1(ii). If $g(i,j)\geq 2a+1$, then, by (7),  $d(v_i)+d(v_j)\geq 2a-1$. Thus, we may assume that $g(i,j)=2a$. Then $f(i,j)\geq 4a$ by (8). The last inequality implies that the arcs $y_ix_j$, $y_jx_i$ are in $D$. Therefore,
$M'_{y,x}:=\{y_ix_j,y_jx_i\} \cup \{y_kx_k\, | \, k\in [1,a]\setminus \{i,j\}\}$ is a new perfect matching from 
$Y$ to $X$ in $D$. Since $x_iy_i\in A(D)$, $x_iy_j\notin A(D)$ and $x_jy_i\notin A(D)$, it follows that 
$s(M'_{y,x})>s(M_{y,x})$, which contradicts the choice of $M_{y,x}$. The claim is proved. \fbox \\\\

We now return to the proof of Lemma 4.3.
 Suppose that two  vertices, say $v_1$ and $v_2$, in  $D^*[M_{y,x}]$ are not adjacent and 
 $$
d(v_1)+d(v_2)\leq 2a-2. \eqno (9) 
 $$
This together with (7), $\overrightarrow{a}[x_1,y_1]\leq 1$ and $\overrightarrow{a}[x_2,y_2]\leq 1$ implies that $ g(1,2)\leq 2a+2$. Therefore, 
$$
2a\leq g(1,2)\leq 2a+2.
$$ 

\textbf{Case 1}. $\overrightarrow{a}[x_1,y_1]=0$.

Then from (7), (9) and the fact that $g(1,2)\geq 2a$, it follows that $\overrightarrow{a}[x_2,y_2]=1$ (i.e., $x_2y_2\in A(D)$) and 
$g(1,2)=2a$. From this and (8) it follows  that 
 $f(1,2)\geq 4a$,
which in turn implies that $y_1x_2\in A(D)$ and $y_2x_1\in A(D)$. These contradicts Claim 1(i) since $x_2y_2\in A(D)$.

\textbf{Case 2}. $\overrightarrow{a}[x_1,y_1]=\overrightarrow{a}[x_2,y_2]=1$, i.e., $x_1y_1\in A(D)$ and $x_2y_2\in A(D)$.

From Claim 1(i) it follows that  $y_1x_2\notin A(D)$ or $y_2x_1\notin A(D)$. If $2a\leq g(1,2)\leq 2a+1$, then from (8) it follows that 
$f(1,2)\geq 4a-1$, which in turn implies that $y_1x_2\in A(D)$ and $y_2x_1\in A(D)$, which is a contradiction. We may therefore assume that $g(1,2)=2a+2$. This and  (8) imply that
$f(1,2)\geq 4a-2$.
 Then, since
$y_1x_2\notin A(D)$ or $y_2x_1\notin A(D)$, it follows that 
 $y_1x_2\in A(D)$ or $y_2x_1\in A(D)$. Without loss of generality, we may assume that $y_1x_2\notin A(D)$ and $y_2x_1\in A(D)$.
Then  $f(1,2)=4a-2$, which in turn implies that  
$$
d^-(x_1)=d^+(y_2)=a \quad \hbox{and} \quad d^-(x_2)=d^+(y_1)=a-1. 
$$
  Therefore, 
$$
y_2\rightarrow \{x_1,x_2,\ldots , x_a\}; \, \{y_1,y_2,\ldots , y_a\}\rightarrow x_1; \,
y_1\rightarrow \{x_1,x_3,x_4,\ldots , x_a\}; \, \{y_2,y_3,\ldots , y_a\}\rightarrow x_2. \eqno (10)
$$
since $y_1x_2\notin A(D)$. Using (10), it is easy to see that for all $i\in [3,a]$, 
$$
M^i_{y,x}=\{y_2x_1,y_ix_2, y_1x_i\}\cup \{y_kx_k \, |\, k\in [3,a]\setminus \{i\}\}
$$
 is a perfect matching from $Y$ to $X$ in $D$. Using the facts that the arcs $x_1y_1$, $x_2y_2$ are in $D$, it is not difficult to see that if for some $i\in [3,a]$, either $x_2y_i\notin A(D)$ or $x_iy_1\notin A(D)$ or  $x_iy_i\in A(D)$, then  $s(M^i_{y,x})>s(M_{y,x})$, which contradicts the choice of $M_{y,x}$.
We may therefore assume that $x_iy_i\notin A(D)$ for all $i\in [3,a]$, and
$$
x_2\rightarrow \{y_2,y_3,\ldots , y_a\} \quad \hbox{and} \quad \{x_3,x_4,\ldots , x_a\}\rightarrow y_1.   
$$
This together with (10) give
$$
x_2\leftrightarrow\{y_2,y_3,\ldots , y_a\} \quad \hbox{and} \quad y_1\leftrightarrow \{x_1,x_3,x_4,\ldots , x_a\}.   \eqno(11) 
$$
Since the vertices $y_1$, $x_2$ are not adjacent, from (11) and Lemma 3.4(i)  it follows that
$$
d^-(y_1)=d^+(x_2)=a-1, \quad d^-(v_1)=d^+(v_2)=a-2.   \eqno (12)
$$

From $g(1,2)=2a+2$, $x_1y_1\in A(D)$ and $x_2y_2\in A(D)$ it follows that
 $$
d(v_1)+d(v_2)= g(1,2)
-4=2a-2.$$
This together with the strongly connectedness of  $D^*[M_{y,x}]$ and  (12)  implies that $d^+(v_1)=d^-(v_2)=1$. This means that $d^+(x_1)=d^-(y_2)=2$. Therefore, $d(x_1)=d(y_2)=a+2$ by (10).

Now for every $i\in [3,a]$ consider the perfect matching $M^i_{y,x}$ and its corresponding digraph $D^*[M^i_{y,x}]$. Notice that 
$s(M_{y,x})=s(M^i_{y,x})=a-2$, the vertices $y_1$, $x_2$ are not adjacent and the arcs $x_iy_i$, $x_1y_2$ are not in $A(D)$. Hence, the vertices $v_1^i=\{y_1,x_i\}$, $v_2^i=\{y_i,x_2\}$ in $D^*[M^i_{y,x}]$ are not adjacent. From Claim 1(ii) it follows that in 
 $D^*[M^i_{y,x}]$ the degree sum of every pair of two distinct non-adjacent vertices, other than $\{v_1^i, v_2^i\}$, at least 
$2a-1$. If in $D^*[M^i_{y,x}]$, $d(v_1^i)+d(v_2^i)\leq 2a-2$, then
by an argument to that in the proof of $d(x_1)=d(y_2)=a+2$, we deduce that $d(x_i)=d(y_i)=a+2$ for all $i\in [3,a]$. Therefore for all $i\in [3,a]$, $3a\leq d(x_1)+d(x_i)\leq 2a+4$, i.e., $a\leq 4$. 

Let $a=3$. Then $x_3y_2\in A(D)$ since $x_1y_2\notin A(D)$ and $d^-(y_1)=2$. Now
using (10) and (11), it is easy to check that 
$x_3y_2x_2y_3x_1y_1x_3$ is a cycle of length 6 in $D$.

Let now $a=4$. By Lemma 4.1, we need  to show that $D$ contains cycles of lengths 6 and 8. From $d(x_4)=6$ and $x_4y_4\notin A(D)$ it follows that $x_4y_2\in A(D)$ or $x_4y_3\in A(D)$. 

Assume that $x_3y_4\in A(D)$. Then using (10) and (11) it is not difficult to see that  
$x_3y_4x_2y_2x_1y_1x_3$ is a cycle of length 6,  and $x_3y_4x_4y_2x_2y_3x_1y_1x_3$ ($x_3y_4x_4y_3x_2y_2x_1y_1x_3$) is a cycle of length  8, when $x_4y_2\in A(D)$ (when $x_4y_3\in A(D)$). 

Assume now that
$x_3y_4\notin A(D)$. Then from  $x_4y_4\notin A(D)$ and $d(y_4)=6$ it follows that 
 $x_1y_4\in A(D)$. Now again using (10) and (11), we see that $x_1y_4x_2y_2x_3y_1x_1$ is a cycle of length 6, 
and $x_1y_4x_4y_2x_2y_3x_3y_1x_1$ ($x_1y_4x_4y_3x_2y_2x_3y_1x_1$)  is a  cycle  length 8, when $x_4y_2\in A(D)$ (when $x_4y_3\in A(D)$). Thus, we have shown that if $a=3$ or $a=4$, then $D$ contains cycles of all lengths $2,4,\ldots ,2a$. This completes the proof of Lemma 4.3. \fbox \\\\

We now ready to complete  the proof of Theorem 1.14.\\

\textbf{Proof of Theorem 1.14.}\\

Let $D$ be a digraph satisfying the conditions of Theorem 1.14. By Lemma 4.1, $D$ contains cycles of lengths 2 and 4. By Lemma 3.1, $D$ contains a perfect matching from $Y$ to $X$. Let $M_{y,x}= \{y_ix_i\in A(D)\, |\, i=1, 2, \ldots , a\}$ be a perfect matching from $Y$ to $X$ in $D$ with the maximum size among the sizes of all the perfect matching from $Y$ to $X$ in $D$.  
   By Lemma 4.2, the digraph $D^*[M_{y,x}]$ either 
 contains cycles of all lengths $2, 4,\ldots , 2a$ or  is strongly connected. Assume that $D^*[M_{y,x}]$ is strongly connected. By Lemma 4.3,
either (i) $a\leq 4$ in which case  $D$ contains cycles of all lengths $2, 4,\ldots , 2a$ or (ii) $d(u)+d(v)\geq 2a-1$ for every pair of non-adjacent vertices $u$, $v$ in $D^*[M_{y,x}]$. Assume that the second case holds. Therefore, by Theorem 1.3, either
(a) $D^*[M_{y,x}]$ contains cycles of every length $k$, $k\in [3,a]$ or (b) $a$ is even and $D^*[M_{y,x}]$ is isomorphic to one of 
$K^*_{a/2,a/2}$, $K^*_{a/2,a/2}\setminus \{e\}$
 or (c) $D^*[M_{y,x}]\in \Phi_a^m$, where $(a+1)/2 < m\leq a-1$.\\

(a). In this case, by Lemma 4.1 and Lemma 3.4(ii), $D$ contains cycles of every length $2k$, $k\in [1,a]$.\\

(b). $D^*[M_{y,x}]$ is isomorphic to 
$K^*_{a/2,a/2}$ or $K^*_{a/2,a/2}\setminus \{e\}$
 with partite sets $\{v_1,v_2,\ldots , v_{a/2}\}$ and
 $\{v_{a/2+1},v_{a/2+2},\ldots , v_{a}\}$. Notice that $a\geq 4$ and  $D^*[M_{y,x}]$ contains cycles of every length $2k$, $k\in [1,a/2]$. Therefore, by Lemma 3.4(ii), $D$ contains cycles of every length $4k$. It remains to show that for any $k\in [1,a/2-1]$, $D$ also contains a cycle of length $4k+2$.

We claim that there exist $p\in [1,a/2]$ and $q\in [a/2+1,a]$ such that $y_px_q\in A(D)$. Assume that this is not the case,
i.e., there is no arc from    
$\{y_1,y_2,\ldots , y_{a/2}\}$ to $\{x_{a/2+1},x_{a/2+2},\ldots , x_{a}\}$. Then, since $D^*[M_{y,x}]$ is isomorphic to 
$K^*_{a/2,a/2}$ or $K^*_{a/2,a/2}\setminus \{e\}$, from the definition of $D^*[M_{y,x}]$ it follows that 
$d^+(y_{1})\leq a/2$, $d^+(y_{a/2})\leq a/2$, $d^-(y_{1})\leq a/2+1$ and $d^-(y_{a/2})\leq a/2+1$. Combining these inequalities, we obtain that  $d(y_1)+d(y_{a/2})\leq 2a+2$, which contradicts condition $A_0$ since $a\geq 4$.

Assume that $D^*[M_{y,x}]$  is isomorphic to $K^*_{a/2,a/2}\setminus \{e\}$.  
Without loss of generality, we may assume that $e=v_av_{a/2}$. From the definition of $D^*[M_{y,x}]$
it follows that  
$\{x_1,x_2,\ldots , x_{a/2}\}\rightarrow \{y_{a/2+1},y_{a/2+2},\ldots , y_{a}\}$ and $D$ contains all possible arcs from 
$\{x_{a/2+1},x_{a/2+2},$ $\ldots , x_{a}\}$ to $\{y_{1},y_{2},\ldots , y_{a/2}\}$ except $x_{a}y_{a/2}$.

 If $p=a/2$ and $q=a$ (i.e., $y_{a/2}x_a\in A(D)$), then 
$y_{1}x_1y_{a/2+1}x_{a/2+1}y_2x_2y_{a/2+2}x_{a/2+2}\ldots y_kx_ky_{a/2+k}$ $x_{a/2+k}y_{a/2}x_ay_{1}$ is a cycle of length $4k+2$, where $k\in [1,a/2-1]$. 
Thus, we may assume that $y_{a/2}x_a\notin A(D)$. Then the vertices $x_a$, $y_{a/2}$ are not adjacent since $x_ay_{a/2}\notin A(D)$.  This together with $d^-(y_{a/2},$ $\{x_1,x_2,\ldots , x_{a/2}\})\leq 1$ give
$d(y_{a/2})\leq 3a/2-1$. 
This together with condition $A_0$ implies that  $d(y_{a/2-1})\geq 3a/2+1$. Therefore, $y_{a/2-1}x_a\in A(D)$ since $d^-(y_{a/2-1},\{x_1,x_2,\ldots , x_{a/2}\})\leq 1$. Now it is not difficult to check that 
$y_{1}x_1y_{a/2+1}x_{a/2+1}y_2x_2y_{a/2+2}x_{a/2+2}$ $\ldots y_kx_ky_{a/2+k}x_{a/2+k}y_{a/2-1}x_ay_{1}$ is a cycle of length $4k+2$ when $k\in [1,a/2-2]$, and 
$y_{1}x_1y_{a/2+1}x_{a/2+1}y_2x_2$ $y_{a/2+2}x_{a/2+2}\ldots x_{a/2-2}y_{a-2}x_{a-2}y_{a/2}x_{a/2}y_{a-1}x_{a-1}y_{a/2-1}$ $x_ay_1$ is a cycle of length $2a-2$ when $a\geq 6$. If $a=4$, then $y_2x_2y_4x_4y_1x_3y_2$ is a cycle of length $6=2a-2$. 

Assume now that $D^*[M_{y,x}]$  is isomorphic to $K^*_{a/2,a/2}$. In this case,  without loss of generality, we assume that $p=a/2$ and $q=a$. By repeating the above argument, when $p=a/2$ and $q=a$, we may show that $D$ contains cycles of all lengths $4k+2$, $k\in [1,a/2-1]$.
\\

(c). $D^*[M_{y,x}]\in \Phi_a^m$. Since $D$ contains cycles of lengths 2, 4 (Lemma 4.1) and every digraph in $\Phi_a^m$ is Hamiltonian, we can assume that $a\geq 4$. Let $V(D^*[M_{y,x}])=\{v_1,v_2,\ldots , v_a\}$ and $v_av_{a-1}\ldots v_2v_1v_a$ is a Hamiltonian cycle in  $D^*[M_{y,x}]$. Therefore, by the definition of $D^*[M_{y,x}]$, for all $i\in [2,a]$, $x_iy_{i-1}\in A(D)$ and $x_1y_a\in A(D)$. From the definition of $\Phi_a^m$
we have $d^+(v_a)=1$ and $d^+(v_{a-1})\leq 2$. These means that $d^+(x_a)\leq 2$ and $d^+(x_{a-1})\leq 3$. These together with $d^-(x_a)\leq a$, $d^-(x_{a-1})\leq a$  and condition $A_0$ imply that
 $$
d(x_a)\leq a+2, \quad  d(x_{a-1})\leq a+3\quad \hbox{and} \quad   3a\leq d(x_a)+d(x_{a-1})\leq 2a+5.  \eqno (13)
$$ 
This implies that $a\leq 5$, i.e., $a=4$ or $a=5$. 

Let $a=5$. Then from (13) it follows that  $d(x_a)+d(x_{a-1})= 2a+5$,
$d^-(x_a)=d^-(x_{a-1})=a$, i.e., $\{y_1,y_2,\ldots , y_a\}\rightarrow \{x_a,x_{a-1}\}$. Therefore, 
$y_2x_5y_4x_4y_3x_3y_2$ (respectively, $y_1x_5y_4x_4y_3x_3y_2x_2y_1$) is a cycle of length 6 (respectively, of length 8).
 
Let $a=4$. In this case we need to show that $D$ contains a cycle of length 6. If  $x_1y_3\in A(D)$ (or $y_2x_1\in A(D)$), then 
$x_1y_3x_3y_2x_2y_1x_1$ (respectively, $x_1y_4x_4y_3x_3y_2x_1$) is a cycle of length 6. We may therefore assume that 
$x_1y_3\notin A(D)$ and $y_2x_1\notin A(D)$. Then $d(x_1)=d(x_4)=6$ since $d(x_4)\leq a+2$ and $d(x_1)+d(x_4)\geq 12$.
Therefore,  $d^-(x_4)=4$, which in turn implies that $y_1x_4\in A(D)$. Hence, $y_1x_4y_3x_3y_2x_2y_1$ is a cycle of length 6. Thus, we have shown that if  $D^*[M_{y,x}]\in \Phi_a^m$, then $a=4$ or $a=5$ and $D$ contains cycles of all lengths $2, 4,\ldots ,2a$.
This completes the proof of the theorem. \fbox \\\\

Our proof does not rely on Theorem 1.10, and thus can be seen as an alternate proof of Theorem 1.8. Note that Theorem 1.9 is an immediate consequence of Theorem 1.14, when $a\geq 3$.\\

\textbf{Note added in proof.} Let $D$ be a balanced bipartite digraph of order $2a\geq 6$. Suppose that $D$ is not a directed cycle and $d(x)+d(y)\geq 3a$ for every pair of vertices $x$, $y$ with a common in-neighbor or a common out-neighbor. 
 Using Theorem 1.10, recently Adamus (arXiv:1708.04674v2 [math.CO] 22 Aug 2017) proved that: 

The digraph 
$D$ either contains cycles of each   even length less than or equal to $2a$ or $d(x)+d(y)\geq 3a$ for every two vertices $x$, $y\in V(D)$ from the same partite set of $D$.\\
 
From this and Theorem 1.14 it follows the following theorem by Adamus.\\
 
\textbf{Theorem 1.15.} (Adamus). Let $D$ be a strongly connected balanced bipartite digraph of order $2a\geq 6$. Suppose that $D$ is not a directed cycle and $d(x)+d(y)\geq 3a$ for every pair of vertices $x$, $y$ with a common in-neighbor or a common out-neighbor. Then   
$D$  contains cycles of each   even length less than or equal to $2a$.


\begin{thebibliography}{25}


\bibitem{[1]}  J. Adamus,  A degree sum condition for hamiltonicity in balanced   bipartite digraphs, Graphs and Combinatorics, vol.  \textbf{33}(1) (2017) 43-51.


\bibitem{[2]}  J. Adamus, L. Adamus and A. Yeo,  On the Meyniel  condition for hamiltonicity  in bipartite digraphs,  Discrete Mathematics and Theoretical Computer Science, \textbf{16}(1) (2014) 293-302. 


\bibitem{[3]} J. Bang-Jensen, Y. Guo, A.Yeo,  A new sufficient condition for a digraph to be Hamiltonian,  Discrete Appl. Math.,  \textbf{95} (1999) 61-72.

\bibitem{[4]} J. Bang-Jensen, G. Gutin, Digraphs: Theory,  Algorithms and Applications, Springer, 2001.

\bibitem{[5]} J. Bang-Jensen, G. Gutin, H. Li,  Sufficient conditions for a digraph to be Hamiltonian,  J. Graph Theory, \textbf{22} (2) (1996) 181-187.

 \bibitem{[6]} A. Benhocine, Pancyclism and Meyniel's conditions, Discrete Math., \textbf{58} (1986) 113-120. 


\bibitem{[7]} J.-C. Bermond, C. Thomassen, Cycles in digraphs-A survey, J. Graph Theory \textbf{5} (1981) 1-43.

\bibitem{[8]} S.Kh. Darbinyan,  On pancyclic digraphs, Preprint of the Computing Center of  Akademy of Sciences of Armenia, (1979) 21 p.    


\bibitem{[9]} S.Kh. Darbinyan,   Pancyclicity of digraphs with the Meyniel condition,  Studia Sci. Math. Hungar., \textbf{20} (1-4) (1985) 95-117) (Ph. D. Thesis, Institute Mathematici Akad. Nauk BSSR, Minsk, 1981). (in Russian). 

   
\bibitem{[10]} S.Kh. Darbinyan, On the pancyclicity of digraphs with large semidegrees,  Akad. Nauk Arm. SSR Dokl., \textbf{83} (4) (1986) 99-101 (see also arXiv: 1111.1841v1 [math.CO] 8 Nov 2011). 

\bibitem{[11]} S.Kh. Darbinyan,  Sufficient conditions for  Hamiltonian cycles in bipartite digraphs, arXiv:1604.08733 [math. CO] 29 Apr 2016. 
   
 \bibitem{[12]} S.Kh. Darbinyan,  Sufficient conditions for a balanced bipartite digraph to be even pancyclic, Discrete Appl. Math., Available online 5-JAN-2018, DOI:10.1016/j.dam.2017.12.013. 

 
\bibitem{[13]} A. Ghouila-Houri,  Une condition suffisante d'existence d'un circuit hamiltonien,  CR Acad. Sci. Paris Ser. \textbf{A-B 25} (1960) 495-497.


\bibitem{[14]} G. Gutin,  Cycles and paths in semicomplete multipartite digraphs, theorems and algorithms: a survey, J. Graph Theory \textbf{19} (4) (1995) 481-505.


\bibitem{[15]} R. H\"{a}ggkvist, C. Thomassen,  On pancyclic digraphs,  J. Combin. Theory Ser. \textbf{B  20} (1), (1976) 20-40.

\bibitem{[16]} D. K\"{u}hn and D. Ostus,  A survey on Hamilton cycles in directed graphs, European J. Combin. \textbf{33} (2012) 750-766.


\bibitem{[17]} M. Meyniel,  Une condition suffisante d'existence d'un circuit hamiltonien dans un graphe oriente,  J. Combin. Theory Ser. \textbf{B  14} (1973) 137-147.

\bibitem{[18]} M. Meszka, New sufficient conditions for bipancyclicity of balanced bipartite digraphs, Discrete Math., in print.

\bibitem{[19]} C. Thomassen,  An Ore-type condition implying a digraph to be pancyclic, Discrete Math. \textbf{19}  (1977) 85-92. 

\bibitem{[20]} R. Wang,  A sufficient conditions for a balanced bipartite digraph to be Hamiltonian, Discrete Mathematics and Theoretical Computer Science, vol. \textbf{19}(3) (2017).
 


\end{thebibliography}
\end{document}